\documentclass{article} 
\usepackage[]{amsmath,amsxtra,amssymb,latexsym, amscd,amsthm, amsfonts,longtable, multicol}
\usepackage{hyperref}
\usepackage{graphicx}
\usepackage{caption}
\usepackage{subcaption}
\usepackage[a4paper,left=3.25cm,right=2.5cm,top=2.5cm,bottom=2.25cm]{geometry}
\date{}
\include{acknowledgement}
\include{abstract}

\begin{document}

\title{An Optimal Control Approach for the Treatment of Hepatitis C Patients}

\author{{\sc Anh-Tuan Nguyen{\footnote{Anh-Tuan Nguyen was with the Faculty of Mathematics, Mechanics and Informatics, Hanoi University of Science, Vietnam National University, 334 Nguyen Trai Street, Thanh Xuan District, Hanoi 100020, Vietnam. He is currently with the imec-Vision Lab, Department of Physics, University of Antwerp, Universiteitsplein 1, Wilrijk 2610, Belgium. Email address: \href{mailto:anh-tuan.nguyen@uantwerpen.be}{\tt anh-tuan.nguyen@uantwerpen.be}.}}}, \quad \sc Hien Tran{\footnote{Hien Tran is with the Department of Mathematics, Center for Research in Scientific Computation, Center for Quantitative Sciences in Biomedicine and Math Multimedia Center, North Carolina State University, Raleigh, NC 27695, United States. Email address: \href{mailto:tran@math.ncsu.edu}{\tt tran@math.ncsu.edu}. \newline \newline \textit{Last edited: \today.}}}}
\maketitle 

\noindent
\begin{abstract}
\noindent In this article, the feasibility of using optimal control theory is studied to develop control theoretic methods for personalized treatment of HCV patients. The mathematical model for HCV progression includes compartments for healthy hepatocytes, infected hepatocytes, infectious virions and noninfectious virions. Methodologies are used from optimal control theory to design and synthesize an open-loop control based treatment regimen for HCV dynamics.\\
\\
\textbf{Key words}: HCV dynamical model, partial virologic response (PVR), infected steady states, optimal control theory.
\end{abstract}

\begin{multicols}{2}
\section{Introduction}
Viral diseases are major causes of human morbidity and mortality worldwide. Global health promotion, therefore, requires the establishment of evidence-based interventions at the individual and population levels aiming at the prevention, treatment and control of viral infections. The design, implementation and impact of such actions depend on a broad understanding of not only host-pathogen interactions in vivo, but also of the intertwining environmental, social and cultural cofactors that may contribute to viral disease development. Multidisciplinary research is thus needed for an effective response to be achieved against these diseases.\\
\\
Hepatitis C virus (HCV) infection stands out among human viral infections for its current importance in global public health. HCV is spread primarily by blood-to-blood contact and affects the liver. The infection is often asymptomatic, but chronic infection can lead to scarring of the liver and ultimately to cirrhosis. Up to 20\% of those with cirrhosis will develop liver cancer. The World Health Organization (WHO) estimates that 180 million people worldwide are affected by the hepatitis C virus. It is globally the most common cause of liver cirrhosis, hepatocellular carcinoma, end-stage liver disease and liver transplantation.\\
\\
The standard treatment for hepatitis C is a 48-week program of weekly treatments with pegylated interferon-alpha and ribavirin that is reviewed after 12 weeks and stopped if the patient is not responding. The results of this approach are mixed. Some patients are cured (sustained virology response (SVR)), some appear to be cured as the end of treatment but the virus returns once the treatment has stopped (relapse) while for others the treatment is not effective (partial virology response (PVR)). This highlights the shortcomings of the current approach of the pharmaceutical industry that could be described as product focused. New products are brought to the market and used as a standard treatment for all patients. However, there is a growing realization in the industry that there is a need to move to a patient focused approach in which the treatment is tailored to the individual patient. The area of personalized medicine has, indeed, been the focus of intense current interest and vigorous research efforts. In essence, personalized medicine seeks to make optimal treatment decisions for an individual patient based on all information available for that patient, including not only genetic and genomic characteristics but also all demographic, physiological, and other clinical factors, thereby allowing the right treatment for the right patient.\\
\\
A precondition for getting better diagnostic procedures and follow-up treatment plans for HCV patients is to obtain better knowledge of the pathophysiological states. One approach to achieve this objective is via patient specific mathematical modeling where knowledge of viral kinetic, pharmacokinetic of the drugs, and various aspects of immune responses is integrated to provide a more comprehensive picture of the biology underlying changes in HCV RNA during therapy. However, for the mathematical models to be clinically relevant, they must be developed in close collaboration with medical doctors and be validated with clinical and laboratory data. In addition, analysis can be performed on the patient specific mathematical models to determine which data gives more information about the system dynamics. Such information can be used for designing better diagnostic and treatment techniques, which in turn have potential to be used for analyzing effects of a given treatment for a given patient. Finally, mathematical models can also be used to predict what elements or biological mechanisms should be studied in more detail and what adjustments of behavior or medical treatment can bring forward positive clinical results.\\
\\
The primary objective of HCV therapy is permanent eradication of the virus. Optimal control theory, a mathematical theory for modifying the behavior of nonlinear dynamical through control of system inputs, is a promising approach to suggesting adaptive HCV treatment strategies. In this research, we will study the feasibility of using optimal control theory to develop control theoretic methods for personalized treatment of HCV patients. A potential problem to consider in this area is whether there exists a drug treatment schedule that can sustain a low viral load and a healthy liver while minimizing the amount of drugs used. In particular, it will be illustrated how optimal control methodology can produce a drug dosing strategy and how this treatment strategy has features of Structured Treatment Interruption (STI). In STI type of treatment, the drug cycles from on treatment to off treatment (drug holiday). The STI type treatment has been considered as an alternative type of treatment for human immunodeficiency virus (HIV). For HIV, there are several reasons to consider an STI type of treatment strategy. First, this type of treatment could potentially stimulate the body's own virus specific immune response. An additional argument for treatment interruption is with individuals with drug resistance strains. A treatment interruption has the potential to allow the drug-sensitive strains to out compete the resistance strains \cite{DAVID}.

\section{HCV infection model}
One of the first mathematical model of Hepatitis C viral dynamics was the model by Neumann in 1998 \cite{NEUMANN}. Neumann developed the model to study the antiviral effects of interferon-$\alpha$ (IFN). The mathematical model is given by
\begin{align}\label{eqn:neumann}
\frac{dT}{dt} &= s - dT - \left(1- \eta\right)\beta VT \\ 
\frac{dI}{dt} &= \left(1 - \eta \right)\beta VT - \delta I \\ 
\frac{dV}{dt} &= \left(1-\epsilon\right) pI - cV,
\end{align}
where $T$ and $I$ are respectively the concentrations of healthy and infected hepatocytes, and $V$ is the concentration of viral load. Healthy hepatocytes are produced at rate $s$ and die with death rate constant $d$. Cells become infected with de novo infection rate constant $\beta$ and, once, infected, die with rate constant $\delta$. Hepatitis C virions are produced by infected cells per day and are cleared with clearance rate constant $c$. The possible effects of IFN in this model are to reduce either the production of virions from infected cells by a fraction $\left(1- \epsilon \right)$ or the de novo rate of infection by a fraction $\left(1 - \eta \right)$. Before IFN therapy, $\epsilon = \eta = 0$. Once therapy is initiated, $\epsilon > 0$ or $\eta > 0$ or both. A major contribution of Neumann et al. was their determination that interferon acts through the $\epsilon$ mechanism rather than the $\eta$ mechanism, i.e, IFN acts by inhibiting the production of visions \cite{JOEY}. \\
\\
Standard care of HCV has since evolved to include, in addition to IFN, the antiviral ribavirin, which acts by rendering some visions noninfectious. To study the response of hepatitis C to this combined treatment, Snoeck et al.\cite{SNOECK} followed Dixit et al.\cite{DIXIT} to extend the original model of Neumann et al. to include a noninfectious virions variable. The model is given by
\begin{align}\label{eqn:snoeck}
\frac{dT}{dt} &= s + r T \left( 1 - \frac{T+I}{T_{max}} \right) - dT - \beta  V_I T \\ 
\frac{dI}{dt} &= \beta  V_I  T + r  I  \left( 1 - \frac{T+I}{T_{max}}\right) - \delta I \\
\frac{dV_I}{dt} &= \left( 1 -\bar{\rho} \right)\left( 1 - \bar{\epsilon} \right) p  I - c  V_I \\
\frac{dV_{NI}}{dt} &= \bar{\rho} \left(1-\bar{\epsilon} \right)  p I - c V_{NI}.
\end{align}
Here, $V_I$ represents the concentration of infectious visions and $V_{NI}$ describes the concentration of noninfectious virions. Another difference is the addition of terms with coefficient $r$ to model the liver's regenerative mechanism. The dynamics of this model after the end of treatment are of interest, so Snoeck et al. include exponential decay of the drug efficacies $\epsilon$ and $\rho$, which correspond to IFN and ribavirin, respectively:
$$\bar{\epsilon} = \epsilon e^{-k(t-t_{end})_+},$$
and
$$\bar{\rho} = \rho e^{-k(t-t_{end})_+},$$
where $t_{end}$ marks the end of treatment, and
$$(a)_+ = \begin{cases}
a &\mbox{if } a\geq 0\\
0 &\mbox{otherwise}
\end{cases}.$$
The steady states of the system were obtained by 
Ransley \cite{MARK} as follows:
\begin{align}\label{eqn:snoeck-steady-states}\nonumber
 {T_u} &= \frac{{{T_{\max }}}}{{2r}}\left( {r - d + \sqrt {{{\left( {r - d} \right)}^2} + \frac{{4rs}}{{{T_{\max }}}}} } \right) \\ \nonumber
 {I_u} &= {V_{Iu}} = {V_{NIu}} = 0 \\ \nonumber
 {T_i} &= \frac{1}{2}\left( { - \frac{{{r^2}D}}{{{A^2}}} + \sqrt {{{\left( {\frac{{{r^2}D}}{{{A^2}}}} \right)}^2} + \frac{{4rs{T_{\max }}}}{{{A^2}}}} } \right) \\ \nonumber
 {I_i} &= {T_i}\left( {\frac{A}{r} - 1} \right) + {T_{\max }} - \frac{{\delta {T_{\max }}}}{r} \\ \nonumber
 {V_{Ii}} &= {\epsilon ^*}{\rho ^*}\frac{p}{c}{I_i} \\ \nonumber
 {V_{NIi}} &= {\epsilon ^*}{\rho}\frac{p}{c}{I_i}, \\ \nonumber
 \end{align}
where $A = \frac{{{\epsilon ^*}{\rho ^*}p\beta {T_{\max }}}}{c}$, ${\rho ^*} = 1 - \rho$, $\epsilon ^* = 1- \epsilon$, $D = \frac{{{T_{\max }}}}{{{r^2}}}\left[ {A\left( {r - \delta } \right) + r\left( {d - \delta } \right)} \right]$. Here, $T_u$, $I_u$, $V_{Iu}$, and $V_{NIu}$ denote the uninfected steady states and  $T_i, I_i, V_{Ii}, V_{NIi}$  represents the infected steady states. \\
\\
The Snoeck et al. model was validated against clinical data \cite{JOEY, SNOECK}. In the following, we will simulate the partial virologic response (PVR) in which the patient responds to the treatment initially and then stops responding while still on treatment. The parameters that correspond to PVR are given in the Figure 1.\\
\\
The initial conditions for our system are not known and are difficult to determine clinically. We followed \cite{MARK} and took the initial conditions to be equal to the infected steady states of the model. The PVR simulation results are depicted in the next 3 figures:
\begin{center}
\includegraphics[width=0.91\linewidth]{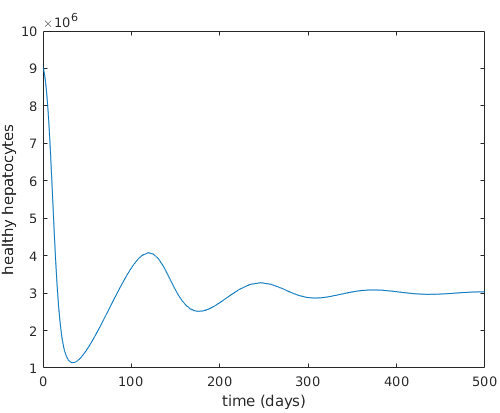}
\end{center}
\begin{center}
\includegraphics[width=0.91\linewidth]{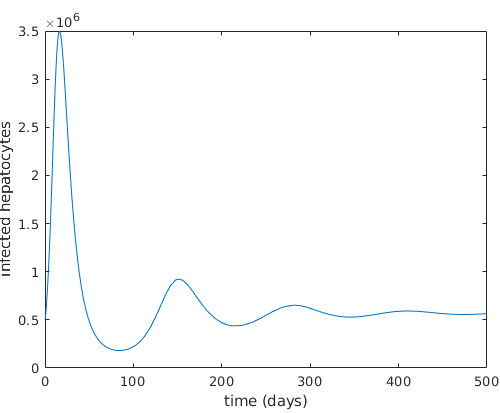}
\end{center}
\begin{center}
\includegraphics[width=0.95\linewidth]{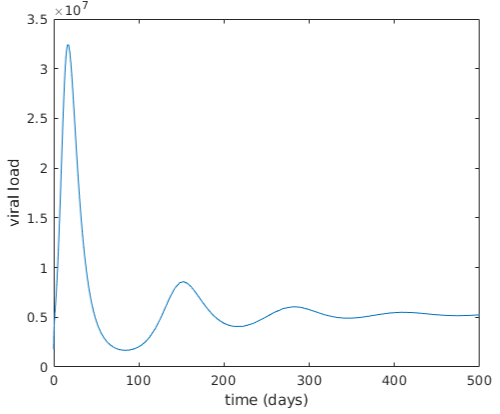} 
\end{center}

\section{An optimal treatment formulation problem}
The problem of finding an effective multi drug therapy might be formulated as an optimal control problem. Let us consider the following optimal control problem:\\
\\
Minimize
$$\int_0^{224}\left[(V_I + V_{NI})^2+I^2-T^2+\epsilon^2+\rho^2\right]\, dt,$$
subject to
\begin{align}\nonumber
\frac{dT}{dt} &= s + r T \left( 1 - \frac{T+I}{T_{max}} \right) - dT - \beta  V_I T \\ \nonumber
\frac{dI}{dt} &= \beta  V_I  T + r  I  \left( 1 - \frac{T+I}{T_{max}}\right) - \delta I \\\nonumber
\frac{dV_I}{dt} &= \left( 1 -\bar{\rho} \right)\left( 1 - \bar{\epsilon} \right) p  I - c  V_I \\\nonumber
\frac{dV_{NI}}{dt} &= \bar{\rho} \left(1-\bar{\epsilon} \right)  pI - c V_{NI},\nonumber
\end{align}
where the components of state should satisfy the initial condition
\begin{equation}\nonumber
     [T(0) \quad I(0) \quad V_I(0)\quad V_{NI}(0)]^T = [T_i \quad I_i \quad V_{Ii} \quad V_{NIi}]
\end{equation}
Here, the viral load $V=V_{I}+V_{NI}$ and the control variables $0\leq \epsilon \leq 1$ and $0\leq \rho \leq 1$ represent the "efficacies" of the drugs IFN and ribavirin, respectively. The goals of the controls are to maximize the healthy hepatocytes while minimizing the viral load, the infected hepatocytes, and the drug usages subject to the HCV dynamics given by the system of equations \eqref{eqn:snoeck}.\\
\\
The treatment is 32 weeks (224 days). The optimal solutions depicted by GPOPS \cite{RAO, SARGENT} in the next five figures showed that the optimal control solutions successfully treated the patient, i.e, sustained virology response (SVR), using a treatment protocol that follows the pattern of structured treatment interruption (that is, on-off-type of treatment). The viral load goes to zero very fast while the concentration of healthy hepatocytes steadily increases and the concentration of infected hepatocytes decreases to zero. This same patient when treated with a continuous dosage ($\epsilon = 6.1382 \times 10^{-1}$, $\rho  = 1.2216 \times 10^{-1}$) of drugs  IFN and ribavirin for 32 weeks only responded to the treatment partially (PVR).\\
\\
\begin{center}
\includegraphics[width=0.91\linewidth]{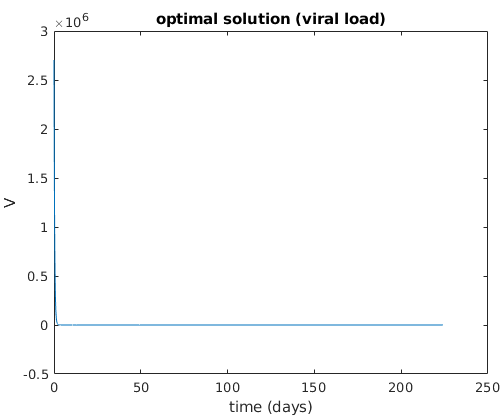}
\end{center}
\begin{center}
\includegraphics[width=0.91\linewidth]{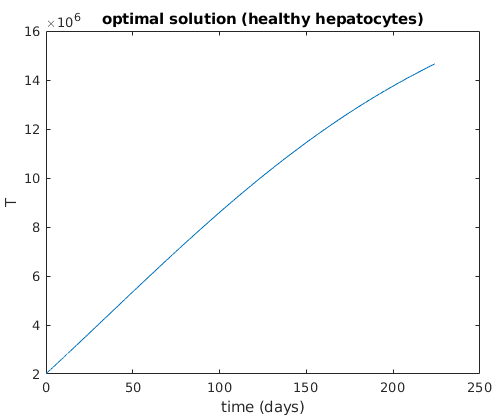}
\end{center}
\begin{center}
\includegraphics[width=0.91\linewidth]{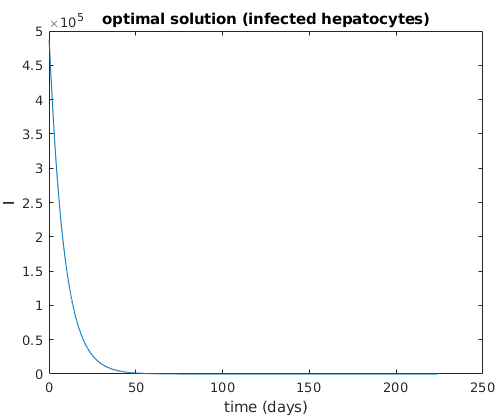}
\end{center}
\begin{center}
\includegraphics[width=0.91\linewidth]{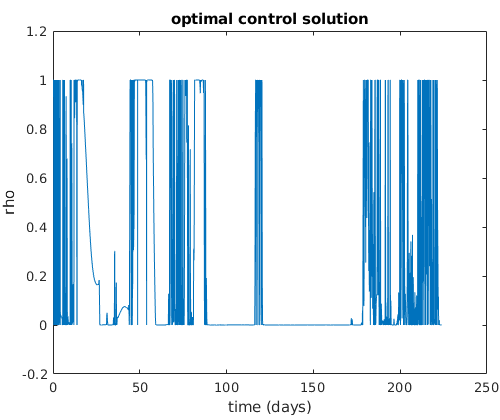} 
\end{center}
\begin{center}
\includegraphics[width=0.91\linewidth]{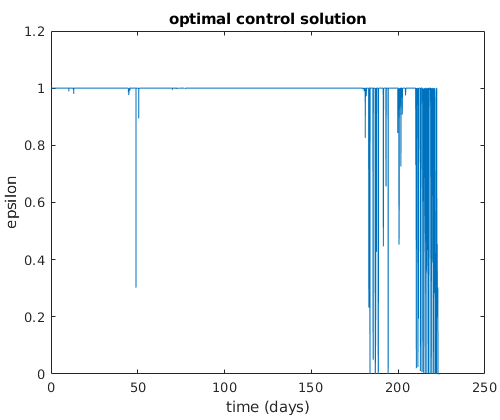}
\end{center}

\noindent
To make sure that the patients remain cured after the treatment has stopped, we simulated the HCV model after the treatment has stopped. The simulated results, which are depicted in the next three figures, showed that viral load and the infected hepatocytes stay below the detected levels while the healthy hepatocytes continued to grow:\\
\\
\begin{center}
\includegraphics[width=0.91\linewidth]{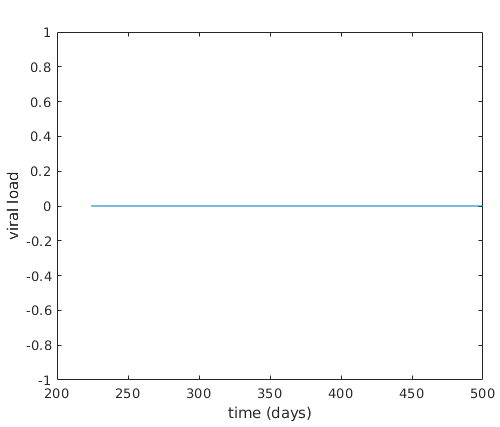}
\end{center}
\begin{center}
\includegraphics[width=0.91\linewidth]{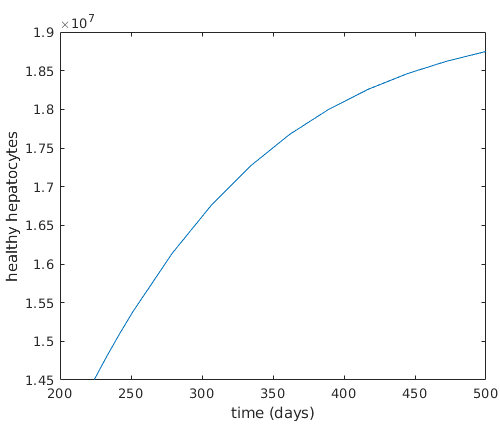}
\end{center}
\begin{center}
\includegraphics[width=0.91\linewidth]{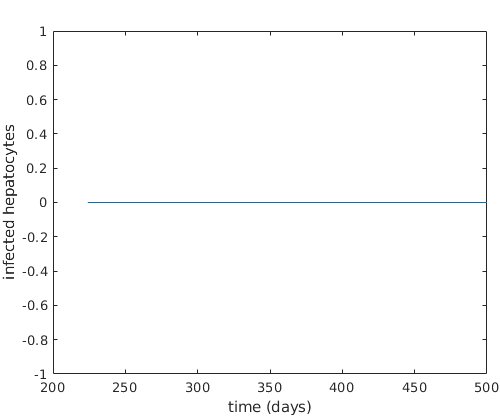}
\end{center}

\section{Conclusion}
We have presented and used methodologies from optimal control theory to design and synthesize an open-loop control based treatment regimen for HCV dynamics. The HCV model is subjected to multiple (IFN and ribavirin) drug treatments as controllers. We demonstrated through numerical simulations that by using a "target tracking" approach, optimal treatment strategies can be designed to successfully cure a patient with HCV (SVR). The optimal treatment strategy followed the pattern of on-off schedule known in the HIV community as structured treatment interruption. It should be emphasized that this same patient when subjected to a continuous treatment regimen only showed partial virologic response (PVR). The potential for an optimal treatment approach for HCV patients is very exciting, but more work is necessary before it can be implemented in clinical settings. For example, more research is needed to incorporate the effects of drug resistance and mutation as well as more accurate mathematical models of HCV dynamics that include the immune response. 

\section*{Acknowledgements}
This research is not financially supported by any research agency or institution. The authors would like to thank Phuong-Thu Trinh for help in MATLAB implementation of validating the method.

\end{multicols}

 \begin{figure}
    \centering
    \begin{longtable}{|c|c|c|c|}
    \hline
    Parameter & Description & Unit & Typical value  \\
    \hline
    $T_{\max}$ & Total number of hepatocytes per ml & H/ml & $18.5\times 10^{+6}$ \\
    s & Hepatocyte production rate & H/ml/day & $61.7 \times 10^{+3}$ \\
    d & Hepatocyte death rate constant & H/day & $0.003$ \\
    r & Hepatocyte proliferation rate constant & H/day & $0.00562$ \\
    p & Virion production rate & IU/day & $25.1$ \\
     k & Antiviral-effect decay constant & day $^{-1}$ & $0.0238$ \\
    $\beta$ & Infection rate constant & H/IU/day & $4.1684\times 10^{-9}$ \\ 
     $\delta$ & Infected cell death rate constant & H/day & $1.2110\times 10^{-1}$\\ 
    c & Virion elimination rate constant & IU/day & $2.7018$ \\
     $\epsilon$ & PEG drug-efficacy & Dimensionless quantity & $6.1382 \times 10^{-1}$ \\
    $\rho$ & RVB drug-efficacy & Dimensionless quantity & $1.2216 \times 10^{-1}$ \\
    \hline
    \end{longtable}
    \vspace{0.5cm}
    \label{fig:parameters}
    \caption{Typical values of the parameters.}
\end{figure}


\begin{thebibliography}{11}
\bibitem{JOEY}{\sc Arthur, J. G., Tran, H. T}. {Feasibility of parameter estimation in hepatitis C viral dynamics models}. \textit{Journal of Inverse and Ill-posed Problems}, \textbf{25}, 1, (2017), pp. {69-80}.

\bibitem{BEELER}{\sc Beeler, S. C., Tran, H. T., Banks, H. T}. {Feedback Control Methodologies for Nonlinear Systems}. \textit{Journal of Optimization, Theory and Applications} \textbf{107}, 1, {(2000)}, {pp. 1-33}.

\bibitem{BANKS}{\sc Banks, H. T., Tran, H. T}. {\it Mathematical and Experimental Modeling for Physical and Biological Processes}. {Chapman \& Hall/CRC Press}, {(2009)}. 

\bibitem{DAVID}{\sc David, J. A}. {\it Optimal Control, Estimation, and Shape Design: Analysis and Applications}. Ph.D Thesis, Department of Mathematics, North Carolina State University, Raleigh, North Carolina, USA. {(2007)}.

\bibitem{DIXIT}{\sc Dixit, N. M., Layden-Almer, J. E., Layden, T. J., Perelson, A. S}. {Modelling how ibavirin improves interferon response rates in Hepatitis C Virus infection}. \textit{Nature}, \textbf{432}, (2004), {pp. 922-924}.

\bibitem{KIRK}{\sc Kirk, D. E}. {\it Optimal Control Theory: An Introduction}. {Prentice-Hall, Inc}. {(1970)}. 

\bibitem{NEUMANN}{\sc Neumann, A. U., Lam, N. P., Dahari, H., Gretch, D. R., Wiley, T. E., Layden, T. J., Perelson, A. S}. {Hepatitis C Viral Dynamics in Vivo and the Antiviral Efficacy of Interferon-$\alpha$ Therapy}. {\textit{Science}, \textbf{282}, 5386}, {(1998)}, {pp. 103-107}.

\bibitem{MARK}{\sc Ransley, M}. {MMath Project: Mathematical Analysis of the HCV Model}. \url{http://www.ucl.ac.uk/\~ucbpran/MMath.pdf} {(2011)}.

\bibitem{RAO}{\sc Rao, A. V., Benson, D., Darby, C. L., Mahon, H., Francolin, C., Patterson, M., Sanders, I., Huntington, G. T}. {\it User's Manual for GPOPS Version 5.0:
A MATLAB Software for Solving Multiple-Phase Optimal Control Problems Using $hp$-Adaptive
Pseudospectral Methods}, (2009).

\bibitem{SARGENT}{\sc Sargent, R. W. H}. {Optimal Control}. {\it Journal of Computational and Applied Mathematics}, \textbf{124}, (2000), pp. 361-371.

\bibitem{SNOECK}{\sc Snoeck, E., Chanu, P., Lavielle, M., Jacqmin, P., Jonsson, E. N., Jorga, K., Goggin, T., Grippo, J., Jumbe, N. L., Frey, N}. {A Comprehensive Hepatitis C Viral Kinetic Model Explaining Cure}. \textit{Clinical Phamacology and Therapeutics} \textbf{87}, 6, (2010), {pp. 706-713}.

\end{thebibliography}
\end{document}